\documentclass[12pt]{article}
\usepackage{epsfig}

\begin{document}
\font \eightrm=cmr8
\font \eightit=cmti8
\renewcommand{\thefootnote}{\fnsymbol{footnote}}
\title{Stochastic Optimal
Prediction with Application to Averaged Euler
Equations\footnotemark[1]}
\author{John Bell\footnotemark[2] \and Alexandre J.\ Chorin\footnotemark[3] \and William
Crutchfield\footnotemark[2]}
\footnotetext[1]{This work was supported in part by
the Office of Science, Office of
Advanced Scientific Computing Research,
Mathematical, Information, and Computational Sciences Division,
Applied Mathematical Sciences Subprogram, of
the U.S.\ Department of Energy,
under Contract No.\ DE-AC03-76SF00098,
and in part by the National Science
Foundation under grant number DMS98-19074.}
\footnotetext[2]
{Lawrence Berkeley National Laboratory, Berkeley, CA  94720.}
\footnotetext[3]
{Mathematics Department, University of California,
Berkeley, CA  94720.}
\date{}
\maketitle

\vskip24pt

\abstract
\vskip14pt

Optimal prediction (OP) methods compensate for a lack of resolution in the
numerical solution of complex problems through the use of an invariant
measure
as a prior measure in the Bayesian sense. In first-order OP, unresolved
information
is approximated by its conditional expectation with respect to the invariant
measure.
In higher-order OP, unresolved information is approximated by a stochastic
estimator,
leading to a system of random or stochastic differential equations.

We explain the ideas through a simple example, and then apply them
to the solution of Averaged Euler equations in two space dimensions.

\newpage

\section{Introduction}

Many problems in mechanics, in particular problems involving turbulence,
cannot be properly resolved because the number of significant degrees of
freedom is too large. The problem of making numerical predictions about the
behavior of systems that have not been
properly resolved has been addressed in \cite{CKK1,CKK2, CKL, HK};
theoretical results can be found in
\cite{hald, HK};
a general introduction to such methods can be found in \cite{chorin7}.
When a system is underresolved,
nothing much can be said without additional information;
in the papers just quoted,
it is assumed that the additional information consists of an invariant
measure on the
space of solutions; this gives rise to an optimal Markovian, deterministic,
approximation.
An invariant measure constitutes additional information because everything
not explicitly known is assumed to be distributed according
to the invariant measure;
it functions like a prior measure in Bayesian statistics \cite{amemiya}.
Unlike what happens
in other areas of application of Bayesian statistics,
nature often provides a rational choice
of prior, invariant measure in the form of a canonical measure.
The use of an invariant measure
gives rise to approximations that are optimal in a sense that we shall specify
below.

In many problems this optimal approximation is still not accurate enough,
and a higher-order,
stochastic, approximation can be derived.
In the present paper we explain these constructions
with the help of a simple example, and then apply them to the solution of
Averaged
Euler equations in two space dimensions.
The main difficulty in higher-order prediction lies
in finding estimators for stochastic processes whose temporal correlations are
determined only empirically.

\section{Optimal and stochastic prediction for Hamiltonian systems}

We present properties associated with general,
even infinite dimensional Hamiltonian systems, in the
simple case of two oscillators,
with position variables $q_1, q_2$ and momenta $p_1,p_2$, and the
Hamiltonian:
\begin{equation}
H=H(q,p)=\frac{1}{2}\left(q_1^2+q_2^2+p_1^2+p_2^2+p_1^2p_2^2\right)
\label{hamiltonian}
\end{equation}
(the ``Hald system").
The equations of the motion of the system are:
\begin{equation}
\matrix{
\frac{dq_1}{dt}&=&p_1+p_1p_2^2,\cr
\frac{dp_1}{dt}&=&-q_1,\cr
\frac{dq_2}{dt}&=&p_2+p_2p_1^2,\cr
\frac{dp_2}{dt}&=&-q_2.}
\label{fullsystem}
\end{equation}

We pretend that 4 equations in 4 unknowns are too difficult to solve on
available computers
but that 2 equations in 2 unknowns are accessible
(A more realistic situation is one where
one has to solve, say, $10^{20}$ equations and one can afford only $10^6$).
Alternately, suppose that for some reason at time $t=0$ we only have values for
$q_1, p_1$
but not for the two other variables.
The question is, how does one write equations for $q_1,p_1$
without computing $q_2,p_2$.

In a standard, Galerkin, approach,
one simply sets all the uncomputed variables to zero;
this results in the system:
$$
\frac{dq_1}{dt}=p_1,\ \  \frac{dp_1}{dt}=-q_1,
\label{galerkin}
$$
which is
not a very good approximation.

Suppose however that although the initial conditions for oscillator 2
are unknown, we do know
that they are drawn from the canonical distribution:
\begin{eqnarray}
P(x_1\leq q_1 <x_1+dx_1, x_2\leq q_2 <x_2+dx_2, \nonumber \\
 y_1\leq p_1<y_1+dy_1, 
y_2\leq p_2<y_2+dy_2)=\nonumber \\
Z^{-1}\exp(-H(x,y)/T)dxdy,
\label{can}
\end{eqnarray}
where $P$ is the probability of the event in parentheses,
$Z$ is a normalization constant that ensures that the sum of all
probabilities is $1$,
$dxdy=dx_1dx_2dy_1dy_2$,
$T$ is a parameter that controls the variance of the samples and is known for
physical reasons as
the temperature,
and $H(x,y)$ is the Hamiltonian function (\ref{hamiltonian}) with $q_1$ replaced
by $x_1$, $p_1$ replaced by $y_1$, etc. Equation (\ref{can})
is often written
in the shorter symbolic form
\begin{equation}
P(q,p)=Z^{-1}\exp(-H(q,p)/T.
\label{canonical}
\end{equation}

One can readily check that this probability distribution is invariant under
the flow
defined by (\ref{fullsystem}), i.e.,
if the initial data are distributed as in (\ref{canonical}), then
the solutions $q_{1,2}, p_{1,2}$ have the same distribution at all later times.
Nature likes this distribution,
and reproduces it often (see any book on statistical mechanics).

Suppose now that the missing initial conditions are drawn from the canonical
distribution
(\ref{canonical}) conditioned by the known information $q_1,p_1$, i.e.,
\begin{equation}
P(x_2\leq q_2<x_2+dx_2, y_2\leq p_2<y_2+dy_2)= 
Z^{-1}\exp(-H_{q_1,p_1}(x_2,y_2))dx_2dy _2,
\label{cond1}
\end{equation}
for some $T$, where $H_{q_1,p_1}$ is $H$ where the values of $x_1,y_1$ have
been given the
fixed, known values of these initial conditions.
In the language of Bayesian estimation \cite{amemiya},
the canonical distribution (\ref{canonical}) is a prior distribution
(what we believe the
distribution to be before we have any data),
and the conditional distribution (\ref{cond1})
is a posterior distribution (the prior distribution modified by what we know
in a special case).
Averages with respect to the conditional distribution (\ref{cond1})
are conditional averages,
and denoted by $E[\cdot|q_1,p_1]$
(the information after the vertical line is what we know,
and the prior distribution is implied). We now approximate equations
(\ref{fullsystem}) by
the optimal prediction (OP) equations:
\begin{equation}
\frac{dq_1}{dt}=E[p_1+p_1p_2^2|q_1,p_1]=p_1+p_1E[p_2^2|q_1,p_1], ~~~
\frac{dp_1}{dt}=-q_1;
\label{op1}
\end{equation}
(clearly $E[p_1|p_1]=p_1$). In our particular case,
\begin{eqnarray*}
E[p_2^2|q_1,p_1]&=&E[p_2^2|p_1]=\cr
&&\int\int
p_2^2\exp\left( -\frac{1}{2}\left(p_1^2+p_2^2+q_1^2+q_2^2+p_1^2p_2^2\right)\right) dq_2
dp_2/\cr
&&\int\int\exp(-(\cdots))dq_2dp_2;
\end{eqnarray*}
(the argument of the exponential is the same in the
denominator as in the numerator).
After obvious cancellations,
$$
E[p^2_2|p_1]=\int p_2^2\exp(-p_2^2/2-p_2^2p_1^2/2)dp_2/\int\exp(\cdots)dp_2
=\frac{1}{1+p_1^2},
$$
(a function of $p_1$).
A general theorem states that in the mean square with respect to the
invariant canonical measure, $E[p^2_2|p_1]$ is the best of approximation
of $p_2^2$ by a function of
$p_1$.
The error in this instance of OP is always smaller than in the Galerkin
approximation above,
though in this instance not by much (see \cite{HK}).
One should think of the system
(\ref{op1}) as producing the average of all solutions obtained by having
initially 
values of $q_1,p_1$ and sampling the other variables from the conditioned
canonical distribution.
An important result due to Hald (\cite{CHK}) states that first-order OP for a
Hamiltonian
system also forms a Hamiltonian system, with a renormalized Hamiltonian
which is minus the logarithm
of the original Hamiltonian averaged over all the ``missing" variables.

First-order OP may be optimal in a mean square sense,
but it may not be good enough in many situations,
and we wish to do better.
In particular, the mean solution of (\ref{fullsystem}) decays, while the
solution of the OP equations (\ref{op1}) does not.

This dichotomy can be understood in several equivalent ways. From 
irreversible statistical mechanics we know that the canonical measure
represents
thermal equilibrium and that the means of all quantities tend to their
equilibrium values
even when they are initially conditioned by partial information;
the symmetry properties of the
Hamiltonian (\ref{hamiltonian}) ensure that the asymptotic mean is zero.
We will now present a second explanation of the decay of the mean
which will motivate our approach to higher-order optimal prediction.
Rewrite the OP equations in the form:
\begin{equation}
\frac{dq_1}{dt}=p_1+p_1z(t), ~~~
\frac{dp_1}{dt}=-q_1,
\label{withz}
\end{equation}
(where $z(t)$ is of course $p_2^2$).
In equations (\ref{op1}) the random function $z(t)$ is
approximated by its conditional mean.
However, in truth $z(t)$ varies from realization to realization of the
initial data,
and we are averaging over systems in which $z(t)$ has a mean value and a
fluctuation around this mean value.
If one thinks of each copy of the system, which conserves energy,
as moving on some constant energy surface,
the surfaces are slightly different for different copies of the system
and the systems move on their surfaces at different rates.
The constant energy surfaces are sphere-like,
the means of these dispersed systems
fall ever closer to the common center of these surfaces,
which is the origin in $qp$ space.
To capture this effect we need to take into account the variability of $z(t)$;
the system (\ref{withz})
with the initial data $q_1,p_1$ can be viewed as a random or stochastic
differential equation;
the problem is that we have yet to figure out what $z(t)$ looks like,
as a function of $t$ and as
a random variable (the randomness coming from the initial conditions).
As we now explain,
the existence of an invariant measure places constraints on $z(t)$ but also
helps in modeling it.
This modeling has to rely on the specific properties of the system under
consideration;
some of the elegant generality of first-order optimal prediction will be lost.
We shall give below an example  of how a term such as $z(t)$ can be estimated.

\section{The Langevin equation and fluctuation / dissipation theorems}

Consider a single particle interacting
with a thermal sea of other particles, the whole being
presumably described in detail by some inaccessibly complicated
Hamiltonian system. 
We wish to describe the evolution of the single particle's velocity without
explicitly calculating the evolution of the sea of particles.
The Langevin equation is a standard approximation of this system:
\begin{equation}
\frac{du}{dt} =-\gamma u+n(t),
\label{langevin1}
\end{equation}
where $u=u(t)$ is the unknown particle velocity, $n(t)$ is white
noise with zero mean
and $\gamma$ is a constant.
The noise $n(t)$ repesents the fluctuating force exerted by the sea of
particles, while the first term on the right hand side represents the
mean force exerted by the sea which opposes the motion of the particle.
The Langevin equation is a 
standard example that shows how an invariant measure constrains the
random forcing
term in a stochastic or random differential equation.
This equation can be solved by elementary means \cite{chandra,chorin7}.
If one thinks of $u$ as the velocity of a particle of mass $1$,
then one should require that asymptotically, as
$t \rightarrow \infty$,
the distribution of $u$ converge to the canonical distribution with density
$Z^{-1}\exp(-u^2/2T)$; this is achieved if
$\gamma=2/T$; this is a ``fluctuation/dissipation'' result
\cite{chandra,chorin7}.
One can understand it as follows: suppose $\gamma=0$; then the variance of $u$
increases with time $t$
(indeed it is proportional to $t$).
If there is no white noise, the initial variability of $u$ decreases
because the term in $\gamma$ is damping.
When the fluctuation/dissipation relation holds, the damping and the
fluctuations imposed by the white noise balance asymptotically so that the invariant distribution is reached.
Application of the fluctuation/dissipation theorems allows the prediction of
the magnitude and direction of the mean force from the
invariant measure and the statistical properties of the random force.
The formula connecting the damping $\gamma$ to the temperature is
based on the assumption that
the autocorrelation of the noise is a delta function, as it is in white noise;
more general fluctuation/dissipation theorems \cite{CHK,mori,zwanzig} will not be used in the present paper.

\section{The Averaged Euler equations}

As an example of the application of these ideas, we consider the Averaged
Euler equations 
\cite{marsden,shkoller1}
in two space dimensions. We choose these equations because an analysis of the
three-dimensional Euler equations
presupposes an extensive discussion of turbulence, while the statistical
mechanics of the
usual two-dimensional Euler equations leads to negative temperature states and
other unusual phenomena.
The two-dimensional Averaged Euler equations 
describe certain temporal
averages of the Euler equations as well as certain viscoelastic flows,
and their statistical mechanics
is compatible with the machinery we have described.

We introduce the operator
$$
A=(1-a^2\Delta)^s,
$$
where $\Delta$ is the Laplace operator,
$a$ is a real constant,
and $s$ is a positive number.
If $s$ is not an integer, $A$ is a pseudo-differential
operator. The Averaged Euler equations are:
\begin{equation}
\frac{\partial A{\bf u}}{\partial
t}+({\bf u}\cdot\nabla)A{\bf u}+ (\nabla {\bf u} )^T \cdot
A{\bf u}=-\nabla p, ~~~
\nabla
\cdot {\bf u}=0,
\label{aeuler}
\end{equation}
where $\bf{u}$ is a vector with components $u_{\alpha}, \alpha=1,2$.
As $a \rightarrow 0$, these equations formally converge to the Euler equations.
We consider a periodic
domain, and expand the $u_{\alpha}$ in Fourier series:
$u_{\alpha}=\sum\hat{u}_{
\alpha,{\bf k}}e^{i{\bf k}\cdot{\bf x}}$,
where ${\bf x}=(x_1,x_2), {\bf k}=(k_1,k_2)$, and $\cdot$ denotes an inner
product.
Substitution into equation (\ref{aeuler}) yields the following equations
of motion for the Fourier
coefficients:
\begin{equation}
\frac{d}{dt} A({\bf k})\hat{u}_{\alpha,{\bf k}}=
-i\sum_{\beta\gamma{\bf k}'} P_{\alpha\beta}({\bf k})
( k'_{\gamma} A({\bf k}) u_{\gamma,{\bf{k-k'}}} u_{\beta,{\bf{k'}}} +
  k'_{\beta} A({\bf k-k'}) u_{\gamma,{\bf{k'}}} u_{\beta,{\bf{k-k'}}} )
\label{foureuler}
\end{equation}
where
$P_{\alpha\beta}({\bf k})=\delta_{\alpha\beta}-\frac{k_{\alpha}k_{\beta}}{{k^2}},$
with 
$\delta$=Kronecker delta, $k^2=k_1^2+k^2_2$,
is the Fourier space projection on the space of divergence-free vectors ($k_1u_{
1,\bf{k}}+k_2u_{2,\bf{k}}=0,$),
$\beta,\gamma$ are component indices and
$A({\bf k})=(1+a^2k^2)^s$ is the Fourier transform of the operator $A$ defined
above.
(This is the straightforward Fourier series of the right-hand side of the
projection form of
equation \ref{aeuler}, with the zero-divergence
condition built-in and the pressure eliminated).

Equation (\ref{aeuler}) conserves an energy, $({\bf u},A{\bf u})$, where the
inner product is the
standard $L_2$ inner product, and an ``enstrophy", $(A\xi,A\xi)$, where $\xi$ is
the vorticity $\xi=\nabla \times {\bf u}$.
Each of these invariants, as well as any of their linear combinations with
positive coefficients, gives rise to an invariant measure with density of
the form $Z^{-1}\exp(-C/T)$, where $C$ is a suitable linear combination.
However,
one can see from general considerations
(\cite{chorin4}) that the energy is irrelevant: if there is no enstrophy in the
expression for the measure
the resulting measure is not ergodic, while if the enstrophy is present the
energy makes little difference.
Thus we consider a measure with density $Z^{-1}\exp(-(A\xi,A\xi)/T)$, in
Fourier variables; (a detailed
example of such a construction is given in \cite{CKL}). The measure is carried 
by divergence-free vectors
(all vectors whose divergence is not zero have probability 0),
and up to the normalizing factor
$Z$ has the density
\begin{equation}
\exp(-\sum k^2(1+a^2 k^2)^{2s}|{\bf\hat{u}}_{\bf k}|^2/T),
\label{measure}
\end{equation}
(For simplicity, we are assuming for the rest of the paper 
that $s=1$ in the operator $A$). This
expression shows why the negative
temperatures of the usual Euler equations do not appear here: they are
necessary in the Euler case to keep
the enstrophy finite \cite{chorin4}; here, for large $k$, 
${\bf\hat{u}}_{\bf k}\sim k^{-3}$, which
makes the enstrophy $\sum k^2|{\bf\hat{u}}|^2$ finite even when the
temperature $T$ is positive.

With this measure, one can see by a simple calculation that if all velocities
are divergence-free,
then
\begin{equation}
E[\hat{u}_{\alpha,{\bf k}}^\ast \hat{u}_{\beta, {\bf k'}}]=
\frac{T\delta_{{\bf k},{\bf k'}} P_{\alpha\beta}({\bf
k})}{k^2(1+a^2k^2)^2}
\ \  \ .
\label{correlations}
\end{equation}
We want to model the evolution of
a small number of Fourier coefficients, those that
satisfy $|{\bf k}|_{\infty}={\rm max}(|k_1|,|k_2|)\le m$,
which we shall call the ``resolved modes''.
We will refer to the remaining modes as ``sampled modes''.
We rewrite the evolution equations for the resolved modes
as sums of
terms that depend only on the resolved modes plus terms that
also involve the sampled modes.
For this purpose we take advantage of the fact that
the right hand side of equation (\ref{foureuler}) is quadratic in
$\bf{u}$ so that its Fourier transform is a convolution of the form
\begin{equation}
\sum\phi_{\alpha\beta}({\bf k},{\bf k'})
\hat{u}_{\alpha,{\bf k'}}\hat{u}_{\beta,{\bf k-k'}},
\label{convolution}
\end{equation}
where the function
$\phi_{\alpha\beta}$ consists of expressions that
guarantee incompressibility and perform the several differentiations; it depends
only on the wave numbers but
not on the amplitudes of the Fourier coefficients.
The terms in the evolution of the resolved modes 
can be divided into three groups:
\begin{enumerate}
\item
Those where both $|{\bf k} - {\bf k'}|_{\infty}\le m$
and $|{\bf k'}|_{\infty}\le m$ (i.e, both
are in the resolved range); we denote their sum by 
$G^1$.
\item
Those where one factor belongs to the resolved range and one does not;
the structure of the convolution
is such that the factor in the sampled range has a wave number such
that $|{\bf k}|_{\infty}\le 2m$.
In other words, resolved modes cannot interact with sampled modes
further away.
If we can model this subset of 
the sampled range, we have all the input we need
to follow the dynamics of the resolved range.
We call their sum 
$G^2$.
\item
Those where neither factor belongs to the resolved range; their sum
is $G^3$.
\end{enumerate}
Thus the Fourier-space evolution equation (\ref{foureuler}) takes the form:
\begin{equation}
\frac{d}{dt}\hat{u}_{\alpha,\bf{k}}=G^1_{\alpha,\bf{k}}+
G^2_{\alpha,\bf{k}} + G^3_{\alpha,\bf{k}}\ \ \ \ \ .
\label{G1G2}
\end{equation}

In first-order optimal prediction $G^2+G^3$ would be approximated by 
its conditional expectation given the
values of the Fourier components
in the resolved range. The Averaged Euler equations share with the
usual Euler equations the remarkable property that when the measure is Gaussian
and based on energy
and/or enstrophy, this conditional expectation is zero provided
${\bf\hat {u}}_ 0=0$.
If ${\bf \hat {u}}_0\ne0$ the conditional expectation takes a simple form,
that we shall not specify
here
because it will not be needed.

\section{Monte Carlo simulations}
\label{sec:monte-carlo}

The system of equations (\ref{foureuler}) is simple enough so that we can
find its mean solution when we have partial data and the remaining data are
drawn from a given distribution, by sampling the 
distribution $N_{ensemble}$ times, 
evolving the system in time, and averaging. 
The system is also 
simple enough that we can determine the statistical properties of the
modes, in particular their time covariances.
The evolution equations are evolved with a
fourth-order Runge-Kutta ODE solver with adaptive time-step control.
A typical simulation requires $10^2-10^3$
runs for reasonable accuracy because
Monte Carlo methods exhibit errors proportional to
$1/\sqrt{N_{ensemble}}$.

Our Monte Carlo simulation lead us to the following observations.
First, the simulation results are only
weakly dependent on the size of the sampled region as long as it
includes $\{ {\bf k} |\  |{\bf k}|_\infty \leq 2m \}$.
The evidence for
this observation will be presented elsewhere.

In the rest of this paper the time correlation function
\begin{equation}
C(\alpha,{\bf k_1},t_1;\beta,{\bf k_2},t_2) = < 
(u_{\alpha,{\bf k_1}}(t_1)-\overline{u_{\alpha,{\bf k_1}}(t_1)})^\ast 
(u_{\beta ,{\bf k_2}}(t_2)-\overline{u_{\beta ,{\bf k_2}}(t_2)}) 
>
\end{equation}
of the sampled modes plays an important role. 
Our numerical simulations lead us to the following observations about
the correlation function:
At time $t=0$, when the probability
distribution is known, the various Fourier components are independent.
We observe that, to a very good approximation, this remains true for $t>0$, and
the time-correlation function is diagonal in $\bf k$.
We may therefore
restrict consideration to the autocorrelation functions of the sampled
modes.  The spatial
index structure of the autocorrelation 
is determined by the requirement that the velocity
fields be divergence free.  We make two further observations about the
correlation functions:
\begin{itemize}
  \item The autocorrelation functions are well approximated in time
  as Gaussians whose width is a function of $\bf k$.
  \item The peak height of the Gaussian is well approximated by the
  magnitude of the correlation function in the invariant measure
  (\ref{correlations}).  Certainly this must be true
  asymptotically as $t\rightarrow\infty$.  
  Our numerical calculations indicate that this is approximately true
  at all times.
\end{itemize}
These observations allow us to describe the 
correlation functions with a single parameter $\sigma({\bf k})$:
\begin{equation}
C(\alpha,{\bf k},t_1;\beta,{\bf p},t_2) \sim
\delta_{\bf k,p} \frac{T}{k^2(1+a^2k^2)^2}
P_{\alpha,\beta} \ e^{(t_1-t_2)^2/\sigma({\bf k})^2}
\label{eq:parametrized-correlation}
\end{equation}

\begin{figure}[tb]
  \epsfysize=3.0truein
  \hfil \epsffile{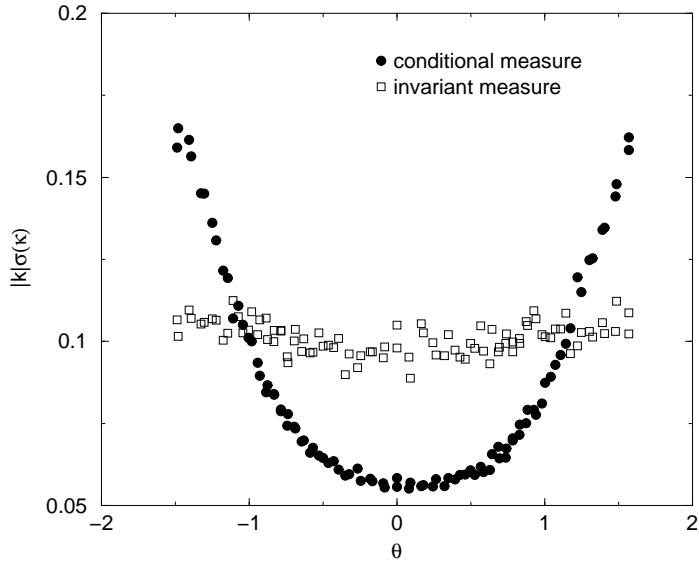} \hfil
  \caption{
  Comparison of autocorrelation widths $\sigma({\bf k})$ in the invariant
  measure and in a specific run with prescribed partial data.  Widths are scaled by
  multiplying by the magnitude of $\bf k$.
  }
  \label{fig:correlation-funct-comparison}
\end{figure}

In general, $\sigma({\bf k})$ will be a function of the direction and
magnitude of
the momentum vector $\bf k$, time $t$,  as well as the initial value
of the resolved
components.  However, in practical application, it will be
advantageous to approximate $\sigma({\bf k})$ by a function of the 
magnitude of $\bf k$ only.  A reasonable approach to this would be
to Monte Carlo simulate the correlation function $C$ with all
$u_{\alpha,{\bf k}}$ chosen from the invariant measure (i.e. no
resolved modes).  This produces
the correlation functions of the invariant measure, which can only be
a function of $|\bf k|$.  Figure
\ref{fig:correlation-funct-comparison} compares $\sigma({\bf k})$ at equilibrium,
i.e., with all the modes sampled, 
to $\sigma({\bf k})$ calculated when the initial measure has specified
resolved modes in the wave-number region $[-5,5] \times [-5,5]$.

Having determined the autocorrelation functions in the prior distribution, 
we shall now use them to approximate specific initial value problems with
prescribed partial data, thus pushing the general methodology of OP to a higher
order in the statistics.  In particular, 
$\sigma({\bf k})$ will henceforth be modeled as a
constant divided by the magnitude of $\bf k$, where the constant is
determined by the invariant measure, as suggested by Figure
\ref{fig:correlation-funct-comparison}.

\section{Approximation of an underresolved system by a stochastic differential equation}

Equation (\ref{G1G2}) shows that the resolved
modes interact with the sampled modes only through the terms $G^2$ and
$G^3$.  
We shall ignore $G^3$ since it is higher order in the sampled modes
and hence suppressed by factors of $1/k^3$.
In the following, we shall use $\hat{u}_{\alpha,\bf k}$ to refer only to the
resolved modes, and $\hat{v}_{\alpha,\bf k}$ to refer to the sampled modes.  The
sampled modes may be written as a mean plus a fluctuation:
$\hat{\bf v}_{\bf k} = \overline{\hat{\bf v}_{\bf k}} + \delta
\hat{\bf v}_{\bf k}$. 
Because $G^2$ is linear in the sampled modes, we may write it
 as $L({\bf \hat{u}})(\overline{\bf \hat{v}}+\delta {\bf \hat{v}})$.  
The evolution equation for the resolved modes may now be written as
\begin{equation}
\frac{d}{dt} \hat{u}_{\alpha,{\bf k}} =
G^1_{\alpha,{\bf k}}(\hat{\bf u}) +
L_{\alpha,{\bf k}; \beta,{\bf p}} \overline{\hat{v}_{\beta,{\bf p}}} +
L_{\alpha,{\bf k}; \beta,{\bf p}} \delta \hat{v}_{\beta,{\bf p}} 
\end{equation}

We will view the 
fluctuating part of the sampled Fourier modes, $\delta \bf \hat{v}$,  
not as a dynamically evolving function of a
random initial configuration, but rather as a random variable with
specific statistical properties, just like $z(t)$ in section 2.
By construction $\delta \bf \hat{v}$ has zero
mean.  The results of the previous section require $\delta \hat{\bf v}_{\bf k}$ to be
uncorrelated with any mode with different $\bf k$, and that its
autocorrelation have the form of 
equation (\ref{eq:parametrized-correlation}).  We
shall not specify any higher order statistics.

The evolution equation above also involves the mean value
$\overline{\bf \hat{v}}$.  We have not presented any empirical observations
about $\overline{\bf \hat{v}}$.  However, by analogy to the Langevin equation,
we can see that the average effect of the $\overline{\bf \hat{v}}$ term in the
evolution equation is to provide a dissipation that counteracts the
fluctuating $\delta \bf \hat{v}$ term.  
We thus write the evolution equation as
\begin{equation}
\frac{d}{dt} \hat{u}_{\alpha,{\bf k}} =
G^1_{\alpha,{\bf k}}(\hat{\bf u}) -
\gamma_{\alpha,{\bf k}; \beta,{\bf p}} \hat{u}_{\beta,{\bf p}}+
L_{\alpha,{\bf k}; \beta,{\bf p}} \delta \hat{v}_{\beta,{\bf p}} 
\label{eq:stochastic-evolution}
\end{equation}
We can use the fluctuation-dissipation
theorem to estimate the dissipation matrix $\gamma$.  Because the
matrix $\gamma$ is dominated by the diagonal elements, in numerical
computations we will
only use the diagonal element:
\begin{equation}
\gamma_{{\bf k},{\bf k}} =
\sum_{\bf k=p+q} 
\frac{T \sqrt{\pi} \sigma({\bf k})
({\bf q}^\perp \cdot  {\bf p})(A({\bf p})p^2-A({\bf q}){ q}^2)^2 |\hat{\bf u}_{\bf p}|^2
}{{q}^4 A({\bf q})^2 p^2 {k}^2 A({\bf k})^2 }
\label{eq:diagonal-gamma}
\end{equation}
where $q^\perp = (q_2,-q_1)$ and q runs over the sampled modes while p
runs over the resolved modes.

The new evolution equation is now a stochastic differential equation
driven by random inputs $\delta \bf v$.  The only Fourier modes that are
evolved by the equations are the resolved modes; the sampled modes are
not evolved.  The statistics of the sampled modes represent the
effect of the unresolved modes.   

The numerical algorithm 
described in section \ref{sec:monte-carlo} 
can be modified for the stochastic differential equation.  
First, only the resolved modes are evolved.
The Fourier amplitudes
$\hat{u}_{\alpha,{\bf k}}$ in the sampled region are chosen
randomly from a population with the correct autocorrelation function
(\ref{eq:parametrized-correlation}).  This correctly models the
$L({\bf \hat{u}})\delta \hat{\bf v}$ term in 
equation (\ref{eq:stochastic-evolution}).
Third, a dissipation term (\ref{eq:diagonal-gamma})
is introduced on the right hand side to
model the $\gamma u$ term in (\ref{eq:stochastic-evolution}); this extracts from
the noise the main part of its mean. 
Unlike what we saw in the Langevin equation of section 3 which has
purely additive noise, we have no
guarantee here that the fluctuation/dissipation formula will take into
account all of the dissipation.
In the evaluation of the dissipation term, we assume that the autocorrelation of each mode
is a constant times a delta function, with  the constant equal to the integral of the
autocorrelation. 
What we have done is tranform a problem with random initial data into a stochastic differential equation. 

An initial test of the stochastic differential equation approach
consists in computing 
a norm of the average solution as a function of time.
We choose as norm the  A-enstrophy of the solution:
\begin{equation}
A_enstrophy({\bf u}) = \sum_{\bf k} k^2 (1+a^2 k^2)^2 |\hat{{\bf u}}_{\bf k}|^2.
\end{equation}
This norm  gives an equal weighting to all the modes. 
The A-enstrophy is, of course, conserved in time in each realisation of the system, but not
on the average ( see the discussion in section 2).

\begin{figure}[tb]
  \epsfysize=3.0truein
  \hfil \epsffile{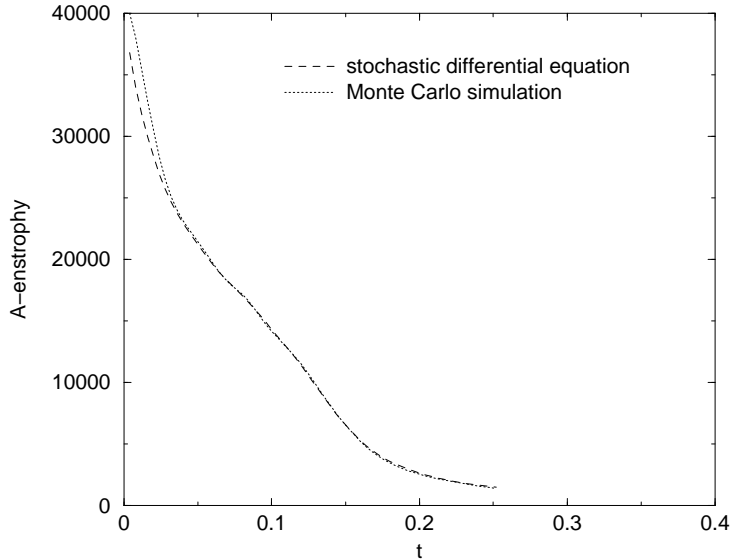} \hfil
  \caption{
  Comparison of decay of the mean A-enstrophy.
  The physical domain is $[0,2\pi]\times[0,2\pi]$.
  The resolved region in wave-number space is $[-5,5]\times[-5,5]$.
  The parameter $a$ in the Average Euler equation was taken to be 1.
  First curve is the true decay as calculated by Monte Carlo.  
  Second curve results from approximating
  system as stochastic differential equation.
  }
  \label{fig:Aenstrophy-decay}
\end{figure}

Figure \ref{fig:Aenstrophy-decay} shows the decay of the A-enstrophy 
of the average solution.
One curve shows the true result obtained from
many Monte Carlo samples of the full equations.  The second curve
results from the stochastic differential equation.
Clearly the decay characteristics are accurately modeled by the
stochastic differential equation.  This is an important improvement
over first order optimal prediction whose evolution is governed by a
renormalized Hamiltonian as noted above, and hence is not dissipative.

\section{Conclusions}

We have used the ideas of optimal prediction to reduce an underresolved 
problem to a stochastic differential equation. This stochastic differential
equation has fewer modes than the full equation, and requires only partial
information about the initial state. On the other hand, it does require
prior knowledge about the statistics of the solution, as we expect in OP 
methods. The full power of higher-order, stochastic OP will appear when we
create effective variance-reduction techniques for the stochastic differential
equation. 
The real test of the ideas will come when we attempt
to solve
Euler and Navier-Stokes equations in three space dimensions; 
the outstanding problem is the
formulation of a reasonable invariant measure in those cases. 
The OP approach
transfers the onus of modeling turbulence from trying to guess relations between
moments to
trying to guess the relevant invariant measures (and awaiting a
mathematical derivation of such measures).
We also expect that the OP machinery will find important uses in other problems
dominated by
complexity, for example in molecular dynamics (see e.g.\cite{K1}).

Note that Figure \ref{fig:Aenstrophy-decay} illustrates a major principle that is often
overlooked: an underresolved conservative system behaves on the average as a 
dissipative system. The importance of this fact for the understanding of turbulence
cannot be overstated.

\section{Acknowledgements}

We would like to thank Prof. G.I. Barenblatt, Prof. O. Hald, Dr. A. Kast,
Prof. R. Kupferman, Prof. D. Levy,
Mr. P. Okunev, and
Prof. B. Turkington, for helpful discussions and for providing much of the
mathematical machinery
upon which this work is based.

\end{document}